\documentclass[12pt]{amsart}
\usepackage{graphics,latexsym,amssymb,amsmath}
\usepackage[dvips]{graphicx}
\usepackage{psfrag,pst-node}
\input{amssym.def}

\newtheorem{thm}{Theorem}
\newtheorem{dfn}{Definition}[section]
\newtheorem{lemma}[dfn]{Lemma}
\newtheorem{rmk}[dfn]{Remark}
\newtheorem{cor}[dfn]{Corollary}
\newtheorem{prop}[dfn]{Proposition}

\newcommand{\Z}{{{\mathbb Z}}}
\newcommand{\N}{{{\mathbb N}}}

\date{}
\begin{document}
\title[$2$-groups with mixed Beauville structures]{An infinite family
  of 2-groups with mixed Beauville structures}

\author[N.~Barker]{Nathan Barker} \author[N.~Boston]{Nigel Boston}
\author[N.~Peyerimhoff]{Norbert Peyerimhoff} \author[A.~Vdovina]{Alina
  Vdovina}

\address[N.~Barker and A.~Vdovina]{School of Mathematics and
  Statistics, Newcastle University, Newcastle-upon-Tyne, NE1 7RU, UK}
\email{nathan.william.barker@googlemail.com}
\email{alina.vdovina@ncl.ac.uk}

\address[N.~Boston]{Department of Mathematics, University of
  Wisconsin, Madison, WI 53706,
  USA} 
\email{boston@math.wisc.edu}

\address[N.~Peyerimhoff]{Department of Mathematical Sciences, Durham
  University, Science Laboratories South Road, Durham, DH1 3LE, UK}
\email{norbert.peyerimhoff@durham.ac.uk}

\begin{abstract}
  We construct an infinite family of triples $(G_k,H_k,T_k)$, where
  $G_k$ are $2$-groups of increasing order, $H_k$ are index-$2$
  subgroups of $G_k$, and $T_k$ are pairs of generators of $H_k$. We
  show that the triples $u_k = (G_k,H_k,T_k)$ are mixed Beauville
  structures if $k$ is not a power of $2$. This is the first known
  infinite family of $2$-groups admitting mixed Beauville
  structures. Moreover, the associated Beauville surface $S(u_3)$ is
  real and, for $k > 3$ not a power of $2$, the Beauville surface
  $S(u_k)$ is not biholomorphic to $\overline{S(u_k)}$.
\end{abstract}

\maketitle

\section{Introduction}

In this article we construct {\em infinitely many} $2$-groups $G_k$
and show that they admit mixed Beauville structures if $k$ is not a
power of $2$.

It was mentioned in \cite{BCG1} that {\em it is rather difficult to
  find a finite group admitting a mixed Beauville structure}. Computer
calculations show that there are no such groups of order $< 2^8$ (see
\cite[Remark 4.2]{BCG2}). By the definition, if a $p$-group admits a
mixed Beauville structure, then $p=2$. Until now, only finitely many
$2$-groups admitting mixed Beauville structures are known. There are
two examples of order $2^8$ in \cite{BCG2}, and five more of orders
$2^{14}, 2^{16}, 2^{19}, 2^{24}, 2^{27}$ in \cite{BBPV}. The family in
this paper is the first known {\em infinite family} of $2$-groups
admitting mixed Beauville structures.

A {\em mixed Beauville structure} of a finite group $G$ is a triple
$(G,H,T)$, where $H$ is an index $2$ subgroup of $G$ and $T =
(h_1,h_2)$ is a pair of elements $h_1,h_2 \in H$ generating $H$ with
particular properties.

Since so little is known about groups admitting mixed Beauville
structures it is generally assumed that they are very rare. Clearly no
simple group can admit a mixed Beauville structure. B. Fairbairn
proved that the same holds true for all almost simple groups $G$
whose derived groups $[G,G]$ are sporadic (see \cite[Theorem 8]{F}). The
only other known construction of groups admitting mixed Beauville
structures was given in \cite{BCG1}. These groups are of the form
$K_{[4]} = (K \times K) \rtimes (\Z/4\Z)$, where $K$ is a group with
particular properties listed in \cite[Lemma 4.5]{BCG1}. The nature of
these other mixed Beauville structures $(K_{[4]},K_{[2]},T=(a,c))$ is
very different from our family of $2$-groups. For example, $\nu(T) =
{\rm ord}(a){\rm ord}(c){\rm ord}(ac)$ contains necessarily two
different primes. Since, for $2$-groups, $\nu(T)$ is necessarily a
power of $2$, this other construction cannot provide examples of {\em
  $2$-groups} admitting mixed Beauville structures.

Our groups $G_k$ are $2$-quotients of a just infinite group $G$ with
seven generators $x_0,\dots,x_6$, acting simply transitively on the
vertices of an $\widetilde A_2$-building. This infinite group first
appeared in \cite{EH}, and then again in \cite{CMSZ} in connection
with buildings. In \cite{PV}, we observed that $G$ has an index $2$
subgroup $H$, generated by $x_0,x_1$, and we used the corresponding index
$2$ quotients $H_k \lhd G_k$ for explicit Cayley graph expander
constructions. The considerations in \cite{PV} showed that $|G_3| =
2^8$ and, for $k \ge 3$,
$$ 
|G_{k+1}| \ge \begin{cases} 8\, |G_k|, & \text{if $k \equiv 0,1 \mod 3$,}\\
  4\, |G_k|, & \text{if $k \equiv 2 \mod 3$.} \end{cases}
$$  
For simplicity of notation, we use the same symbols $x_i$ for the
generators of $G$ and their images in the finite quotients $G_k$.
 
Any mixed Beauville structure $u = (G,H,T)$ gives rise to a Beauville
surface $S(u) \cong (C_T \times C_T)/G$ of mixed type. A natural
question is whether this Beauville surface $S(u)$ is real. An algebraic
surface $S$ is called {\em real} if there is a biholomorphism $\sigma:
S \to \overline{S}$ with $\sigma^2 = {\rm id}$. For the details we
refer, e.g., to the papers \cite{BCG1} and \cite{BCG2}.

Let us now formulate the main result of this paper.

\begin{thm} \label{thm:main} 
  Let $k \ge 3$ be {\em not} a power of $2$ and $T_k = (x_0,x_1) \in
  H_k \times H_k$. Then the triple $u_k=(G_k,H_k,T_k)$ is a mixed
  Beauville structure. Moreover
  \begin{itemize}
  \item[(i)] The mixed Beauville surface $S(u_3)$ is real. 
  \item[(ii)] For every $k > 3$ not a power of $2$, the Beauville
    surface $S(u_k)$ is not biholomorphic to its complex conjugate
    $\overline{S(u_k)}$.
  \end{itemize}
\end{thm} 

For the proof, we realise $G$ as a group of (finite band) upper
triangular infinite Toeplitz matrices. The $2$-quotients $G_k$ are
obtained via truncations of these matrices at their $(k+1)$-th upper
diagonal, and they have a certain nilpotency structure. Our proof
exploits this nilpotency structure as well as subtle periodicity
properties of these matrices. It also becomes transparent via these
periodicity properties why, in the above theorem, $k \ge 3$ must
necessarily avoid the powers of $2$.

Let us explain the difference between the results in \cite{BBPV} and
in this article: In \cite{BBPV}, we used the computational algebra
system Magma to check that the {\em first six groups} of an infinite
family of $2$-groups admit mixed Beauville structures, which led us to
conjecture that this holds true for the full infinite family. In this
paper, we provide a {\em rigorous theoretical proof} that an infinite
family of $2$-groups admit mixed Beauville structures. In view of the
final Remark \ref{rmk:reflection}, it is very surprising that all our
groups (except for $G_{2^j}$ with $j \in \N_0$) admit mixed Beauville
structures. Moreover, there is overwhelming evidence that the families
of groups in both papers agree, and it has been verified
computationally for the first $100$ groups in both families that they
are pairwise isomorphic (see \cite[Conjecture 1]{PV}).

Let us finish our introduction with the following question: {\em For
  which $2$-groups $H$ does there exist a group $G \supset H$ and a
  choice $T \in H \times H$ such that $(G,H,T)$ is a mixed Beauville
  structure?} Both examples of groups of order $2^8$ listed in
\cite[Thm 0.1]{BCG2} and admitting mixed Beauville structures have the
same index $2$-subgroup which agrees with our group $H_3$. The five
other examples in \cite{BBPV} agree with our examples
$H_5,H_6,H_7,H_9,H_{10}$. It would be interesting to know whether
there are any other $2$-groups $H$ giving rise to mixed Beauville
structures $(G,H,T)$, and which do not agree with one of our groups
$H_k$.

\bigskip

{\sc Acknowledgement:} The first author likes to thank Uzi Vishne for
useful correspondences. The research of Nigel Boston is supported by
the NSA Grant MSN115460.

\section{Mixed Beauville structures and associated surfaces}

The following presentation follows \cite{BCG1} closely. Let $G$ be a
finite group and $H \subset G$ be a subgroup of index $2$. For $x \in
H$ let
$$ \Sigma(x) := \{ h x^j h^{-1} \mid h \in H, j \ge 0 \}, $$
i.e., $\Sigma(x)$ is the union of all conjugates of the cyclic
subgroup generated by $x$. For $T = (x_0,x_1) \in H \times H$, we define
$$ \Sigma(T) := \Sigma(x_0) \cup \Sigma(x_1) \cup \Sigma((x_0 x_1)^{-1}). $$ 

A {\em mixed Beauville structure} is a triple $(G,H,T)$ with $T =
(x_0,x_1)$ satisfying the following properties:
\begin{itemize}
  \item[(A)] $x_0,x_1$ generate the group $H$.
  \item[(B)] There exists $g_0 \in G \backslash H$ such that
  $g_0 \Sigma(T) g_0^{-1} \cap \Sigma(T) = \{ {\rm id} \}$.
  \item[(C)] For all $g \in G \backslash H$ we have $g^2 \not\in \Sigma(T)$.
\end{itemize} 

Next, we explain how to construct the Beauville surface $S=S(u)$
associated to a mixed Beauville structure $u=(G,H,T=(x_0,x_1))$. Let
$P_0,P_1,P_2 \in \mathbb{P}^1$ be a sequence of points ordered
counterclockwise around a base point $O \in \mathbb{P}^1$ and, for $0
\le i \le 2$, let $\gamma_i \in \pi_1({\mathbb P}^1 \backslash \{
P_0,P_1,P_2 \},O)$ be represented by a simple counterclockwise loop
around $P_i$ such that $\gamma_0 \gamma_1 \gamma_2 = {\rm id}$. By
Riemann's existence theorem, there exists a surjective homomorphism
$$ \Phi: \pi_1({\mathbb P}^1 \backslash \{ P_0,P_1,P_2 \},O) \to H $$
with $\Phi(\gamma_0) = x_0$ and $\Phi(\gamma_1) = x_1$, and a Galois
covering $\lambda_T: C_T \to {\mathbb P}^1$, ramified only in $\{ P_0,
P_1, P_2 \}$, with ramification indices equal to the orders of the
elements $x_0, x_1, x_0 x_1$. These data induce a well defined
action of $H$ on the curve $C_T$, and by the Riemann-Hurwitz formula,
we have
$$ g(C_T) = 1 + \frac{|H|}{2}\left(1- \frac{1}{{\rm ord}(x_0)} - 
\frac{1}{{\rm ord}(x_1)} - \frac{1}{{\rm ord}(x_0 x_1)}\right). $$ 
Let $\varphi_g: H \to H$ be conjugation with $g$, i.e., $\varphi_g(x)
= g x g^{-1}$. We then define a $G$-action on $C_T \times C_T$ by
$$ x(z_1,z_2) = (x z_1, \varphi_g(x) z_2), \qquad g(z_1,z_2) = (z_2,g^2 z_1), $$
for all $x \in H$ and $(z_1,z_2) \in C_T \times C_T$. This action is
fixed point free, and the quotient $(C_T \times C_T) / G$ is the
associated mixed Beauville surface $S$. By the Theorem of
Zeuthen-Segre, we have for the topological Euler number
\begin{eqnarray*}
e(S) &=& \frac{4 (g(C_T)-1)^2}{|H|} \\
&=& |H| \left(1- \frac{1}{{\rm ord}(x_0)} - \frac{1}{{\rm ord}(x_1)} - 
\frac{1}{{\rm ord}(x_0 x_1)}\right)^2, 
\end{eqnarray*}
as well as the relations (see \cite[Theorem 3.4]{Cat}),
$$ \chi(S) = \frac{e(S)}{4} = \frac{{K_S}^2}{8}, $$
where ${K_S}^2$ is the self-intersection number of the canonical
divisor and $\chi(S) = 1 + p_g(S) - q(S)$ is the holomorphic
Euler-Poincar{\'e} characteristic of $S$.

Let us briefly indicate how we prove the {\em reality statements}
(i),(ii) for the mixed Beauville surfaces in Theorem \ref{thm:main}:
For $T = (c,a) \in H \times H$ let $T^{-1} = (c^{-1},a^{-1})$. Every
mixed Beauville structure $u = (G,H,T)$ gives rise to another mixed
Beauville structure $\iota(u) = (G,H,T^{-1})$, and we have
$S(\iota(u)) = \overline{S(u)}$ (see \cite[(39)]{BCG1}). Let ${\mathbb
  M}(G) = \{ (G,H,(c,a)) \}$ denote the set of all mixed Beauville
structures of $G$. Every automorphism $\psi \in {\rm Aut}(G)$ induces
a map $\sigma_{\psi}$ on ${\mathbb M}(G)$ via
$$ \sigma_{\psi}(G,H,(c,a)) = (G,\psi(H),(\psi(c),\psi(a))). $$
Moreover, in accordance with \cite[(11) and (32)]{BCG1},
let $\sigma_3,\sigma_4$ be maps on ${\mathbb M}(G)$, defined by
\begin{eqnarray*}
\sigma_3(G,H,(c,a)) &=& (G,H,(a,c)), \\
\sigma_4(G,H,(c,a)) &=& (G,H,(c,c^{-1}a^{-1})),
\end{eqnarray*}
and $A_{\mathbb M}(G)$ be the group generated by the maps $\sigma_\psi$ ($\psi
\in {\rm Aut}(G)$) and $\sigma_3,\sigma_4$. Then we have the following
facts (see \cite[Prop. 4.7]{BCG1}):
\begin{itemize}
\item[(a)] $S(u)$ is biholomorphic to $\overline{S(u)}$ iff $\iota(u)
  \in A_{\mathbb M}(G)u$.
\item[(b)] $S(u)$ is real iff $\iota(u) = \rho(u)$ for some $\rho \in
  A_{\mathbb M}(G)$ with $\rho(\iota(u)) = u$.
\end{itemize}
Choosing the mixed Beauville structures $u_k$ from Theorem
\ref{thm:main}, we find an automorphism $\psi: G_3 \to G_3$, uniquely
defined by $\psi(x_0) = x_0^{-1}$, $\psi(x_1) = x_1^{-1}$ and
$\psi(x_2) = x_0^{-1} x_2 x_0$. This implies $\iota(u_3) =
\sigma_\psi(u_3)$ and $\sigma_\psi(\iota(u_3)) = u_3$, and it follows
from (b) that $S(u_3)$ is real. On the other hand, for $k > 3$ and not
a power of $2$, we show that there is no homomorphism $\psi: H_k \to
H_k$ satisfying
\begin{multline*}
(\psi(x_0),\psi(x_1)) \in \{ (x_0^{-1},x_1^{-1}), (x_1 x_0,x_0^{-1}),
(x_1^{-1},x_1 x_0), \\ (x_1^{-1},x_0^{-1}), (x_0^{-1},x_1 x_0),
(x_1 x_0,x_1^{-1}) \}.
\end{multline*}
Using \cite[Lemma 2.4]{BCG1} and the criterion (a) above, this implies
that $S(u_k)$ cannot be biholomorphic to $\overline{S(u_k)}$. (Note that
our pair $(x_0,x_1)$ corresponds, in the notation of \cite{BCG1}, to
the pair $(c,a)$.)

\section{The $2$-groups $G_k$ and $H_k$}

Let $\mathcal K$ be the simplicial complex constructed from the
following $7$ triangles by identifying sides with
the same labels $x_i$. 

\begin{figure}[h]
  \begin{center}      
      \psfrag{x0}{$x_0$}
      \psfrag{x1}{$x_1$}
      \psfrag{x2}{$x_2$}
      \psfrag{x3}{$x_3$}
      \psfrag{x4}{$x_4$}
      \psfrag{x5}{$x_5$}
      \psfrag{x6}{$x_6$}
     \includegraphics[width=\textwidth]{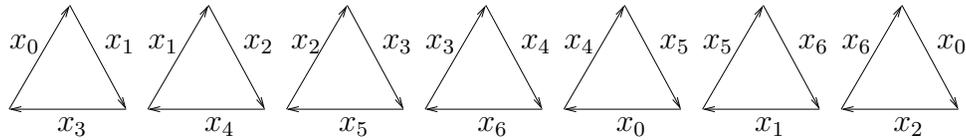}
  \end{center}
   \caption{Labeling scheme for the simplicial complex ${\mathcal K}$}
   \label{labelscheme}
\end{figure} 

It is easily checked that the vertices of all triangles are
identified, and that the fundamental group $\pi_1({\mathcal K})$ is
isomorphic to the infinite abstract group
\begin{equation} \label{eq:gppres}
G = \left\langle x_0, \dots, x_6 \mid \text{$x_i x_{i+1} x_{i+3} = {\rm id}$
    for $i=0,\dots,6$} \right\rangle, 
\end{equation}
where $i,i+1$ and $i+3$ are taken modulo $7$. Realising the triangles
as equilateral Euclidean triangles, we can view the universal covering
of $\mathcal K$ as a thick Euclidean building of type
$\widetilde{A_2}$, on which $G$ acts via covering transformations.

Note that the presentation \eqref{eq:gppres} is a presentation of $G$
by $7$ generators and $7$ relations. It is easy to see that $G$ is
already generated by the three elements $x_0,x_1,x_2$. Let $H \subset
G$ be the subgroup generated by the two elements $x_0,x_1$. Then $H$
is an index-$2$ subgroup of $G$ (see \cite[Prop. 2.1]{PV}). The groups $G_k$ and
$H_k$ will be finite $2$-quotients of these groups $G$ and $H$.

We now recall the {\em faithful representation} of $G$ by infinite upper triangular
matrices given in \cite{PV}, where every element $x \in G$ is represented as
\begin{equation} \label{eq:xipm}
x = \begin{pmatrix} 
1 & a_{11} & a_{21} & \dots & a_{k1} & 0 & 0 & \dots & \dots \\ 
0 & 1 & a_{12} & a_{22} & \dots & a_{k2} & 0 & \ddots \\
0 & 0 & 1 & a_{13} & a_{23} & \dots & a_{k3} & 0 & \ddots \\
\vdots & \ddots & 0 & 1 & a_{11} & a_{21} & \dots & a_{k1} & \ddots \\
\vdots &        & \ddots & \ddots & 1 & a_{12} & a_{22} & \dots & \ddots \\
\vdots &        &        & \ddots & \ddots & 1 & a_{13} & a_{23} & \ddots \\
\vdots &        &        &        & \ddots & \ddots & \ddots & \ddots & \ddots
\end{pmatrix}, 
\end{equation}
and each entry $a_{ij}$ is a matrix in $M(3,{\mathbb F}_2)$ (and $0$
and $1$ stand for the zero and the identity matrix in $M(3,{\mathbb
  F}_2)$). Note that the matrix representation \eqref{eq:xipm} has
only finitely many non-zero upper diagonals. Moreover, the entries in
every diagonal are repeating with period $3$.

Let us now introduce a concise notation for these matrices: The
entries $a_{j1}, a_{j2}, a_{j3} \in M(3,{\mathbb F}_2)$ of the $j$-th
upper diagonal in the matrix representation can be combined to a $3
\times 9$ matrix, which we denote by $a_j = [a_{j1}, a_{j2},
a_{j3}]$. (Conversely, we refer to the three $3 \times 3$ matrices
constituting a $3 \times 9$ matrix $a_j$ by $a_j(1),a_j(2),a_j(3) \in
M(3,{\mathbb F}_2)$.) We can then write the matrix in \eqref{eq:xipm} as
$$ M_0(a_1,\dots,a_k) = M_0([a_{11},a_{12},a_{13}],\dots,[a_{k1},a_{k2},a_{k3}]). $$
If the first $l \ge 1$ upper diagonals of a matrix
$M_0(a_1,\dots,a_k)$ are zero, we use also the notation
$M_l(a_{l+1},\dots,a_k)$. Since the presentation with $3 \times 9$
matrices is still not very concise, we translate every matrix $a =
(u_{ij}) \in M(3,{\mathbb F}_2)$ into the non-negative integer
$$ A = 256 u_{11} + 128 u_{12} + 64 u_{13} + 32 u_{21} + 16 u_{22} + 8 u_{23} + 4 u_{31} + 
2 u_{32} + u_{33}, $$ 
and represent the $3 \times 9$ matrix $[a_{j1}, a_{j2}, a_{j3}]$ by
the triple $A_j = [A_{j1},A_{j2},A_{j3}]$ with $0 \le A_1, A_2, A_3 \le
511$. Therefore, another way to write the matrix in \eqref{eq:xipm} is
$$ M_0(A_1,\dots,A_k) = M_0([A_{11},A_{12},A_{13}],\dots,[A_{k1},A_{k2},A_{k3}]). $$
The matrices corresponding to the generators
$x_0,x_1,x_2$ are in this notation:
\begin{eqnarray*}
x_0 &=& M_0([11,11,11],[17,17,17],[26,26,26],[11,11,0],[17,0,0]), \\
x_1 &=& M_0([23,224,138],[59,136,495],[26,488,227],[23,224,0],[59,0,0]),\\
x_2 &=& M_0([46,68,217],[12,194,363],[26,326,77],[46,68,0],[12,0,0]).
\end{eqnarray*}

The proofs of the explicit formulas in the following Lemma are
straightforward (see \cite{PV}). Note that (b) is a refinement of
\cite[Prop. 2.5]{PV}. These formulas are crucial for our later
considerations:

\begin{lemma} \label{lem:crucform}
  Note that in the following formulas all entries $j$ in $a_{1/2}(j)$, $b_1(j)$,
  $c_{1/2}(j)$ are taken mod $3$ and chosen to be in the range $\{1,2,3\}$.
  \begin{itemize}
  \item[(a)] Let $k, j \ge 0$ and $M_1 = M_k(a_1,a_2,\dots)$ and $M_2
    = M_{k+j}(b_1,\dots)$. Then both products $M_1 M_2$ and $M_2 M_1$
    are of the form
    $$ M_k(a_1,a_2,\dots,a_{j-1},a_j + b_1,\dots). $$
  \item[(b)] We have
    $$ M_k(a_1,a_2,\dots)^2 = M_{2k+1}(c_1,c_2,\dots), $$
    with $c_1(i) = a_1(i) a_1(k+i+1)$ and
    $$ c_2(i) = a_1(i) a_2(k+i+1) + a_2(i) a_1(k+i+2). $$
  \item[(c)] We have
    $$ M_0(b_1,\dots)^{-1} M_k(a_1,a_2,\dots) M_0(b_1,\dots) = M_k(a_1,c_2,\dots) $$
    with $c_2(i) = a_2(i) + b_1(i) a_1(i+1) + a_1(i) b_1(k+i+1)$.
  \end{itemize}
\end{lemma}

Let $G^k$ and $H^k$ be the subgroups of all elements in $G$ and $H$
with vanishing first $k$ upper diagonals (i.e., these elements are of
the form $M_k(a_1,\dots)$). Then $G^k$ and $H^k$ are normal subgroups of
$G$ and $H$, and our groups $G_k$ and $H_k$ are the quotients $G / G^k$ and
$H / H^k$. We can think of $G_k$ and $H_k$ as truncations of the matrix groups
$G$ and $H$ at their $(k+1)$-st upper diagonal. The finiteness of these quotients
follows then easily from the $3$-periodicity of the diagonals. 

\begin{rmk} Another way to generate quotients of $G$ (and $H$) is via
  the lower-exponent-$2$ series
  $$ G = \lambda_0(G) \supset \lambda_1(G) \supset \dots $$
  where $\lambda_{i+1}(G) = [\lambda_i(G),G](\lambda_i(G))^2$. The
  quotients $G/\lambda_k(G)$ are finite $2$-groups. It follows from
  \cite[Prop 2.5]{PV} that $\lambda_k(G) \subset G^k$. Magma computations
  show for all indices $k \le 100$ (see \cite{PV}) that $\lambda_k(G) \cong G^k$ and
  $$ \log_2 [ \gamma_k(G):\gamma_{k+1}(G) ] = \begin{cases} 
    3, \text{if $k \equiv 0,1 \mod 3$}, \\ 2, \text{if $k \equiv 2 \mod 3$}. 
    \end{cases} $$ 
  We conjecture (see \cite[Conjectures 1 and 2]{PV}) that these facts hold
  true for all $k$, which would mean that the group $G$ {\em has
  finite width $3$} (see \cite{KLGP} for definitions).
\end{rmk}

\section{Powers of the generators}

This and all the following sections are dedicated to the proof that
the triple $(G_k,H_k,T_k)$ satisfies the conditions (A), (B) and (C)
of a mixed Beauville structure if $k$ is not a power of $2$.  The
triples
$$ [x_0,x_1,x=(x_0 x_1)^{-1}] \,\, \text{and} \,\, 
[y_0=x_2 x_0 x_2^{-1},y_1=x_2 x_1 x_2^{-1},y=(y_0 y_1)^{-1}] $$ are
both {\em spherical systems of generators} of the group $H$ (see,
e.g., \cite{BCG2} for this notion). A crucial step towards the proof
of Theorem \ref{thm:main} is the explicit determination of the first
two non-trivial diagonals of {\em all powers} of each of the elements
$x_0,x_1,x,y_0,y_1,y$. By the {\em first two non-trivial diagonals} of
a matrix $M_0(a_1,a_2,\dots) \neq {\rm id}$ we mean the pair $a_k,
a_{k+1}$ with $a_1=\dots=a_{k-1}=0$ and $a_k \neq 0$. Moreover, we
call $a_k$ the {\em leading diagonal} of this matrix. In fact, it
turns out that -- in all considerations of this paper -- only a good
understanding of the first two non-trivial diagonals is needed and
that the higher diagonals can be ignored.

Let us focus on the powers of the elements
\begin{eqnarray*}
x = (x_0 x_1)^{-1} &=& M_0([28,235,129],[29,211,263],\dots), \\
y = x_2 x x_2^{-1} &=& M_0([28,235,129],[58,3,445],\dots)
\end{eqnarray*} 
for reasons of illustration (the analogous results for the powers of
the pairs $(x_0,y_0=x_2 x_0 x_2^{-1})$ and $(x_1,y_1=x_2 x_1
x_2^{-1})$ will be given at the end of this section).

It is remarkable that the first two non-trivial diagonals of the
$2$-powers of $x$ and $y$ repeat with a {\em periodicity of $2$}. This
is the content of the following proposition and can be verified by a
straightforward calculation using Lemma \ref{lem:crucform}(b):

\begin{prop} \label{prop:2power}
  We have for all $j \ge 0$:
  \begin{eqnarray*}
  x^{2^{2j+1}} &=& M_{2^{2j+1}-1}([51,89,196],[0,0,0],\dots), \\
  y^{2^{2j+1}} &=& M_{2^{2j+1}-1}([51,89,196],[0,157,106],\dots), \\
  x^{2^{2j+2}} &=& M_{2^{2j}+1}([28,235,129],[0,0,0],\dots), \\
  y^{2^{2j+2}} &=& M_{2^{2j}+1}([28,235,129],[39,208,186],\dots).
  \end{eqnarray*}
\end{prop}

\begin{rmk} The group $G$ has more remarkable properties. In
  \cite[Prop. 2.6]{PV}, we present a certain $3$-periodicity of
  commutators. Another interesting property is that the subgroup
  generated by the squares $x_0^2,x_1^2,\dots,x_6^2$ of the seven
  generators is isomorphic to $G$ (see \cite[p. 308]{EH}).
\end{rmk}

The next remark explains why the statement in Theorem \ref{thm:main}
cannot hold for powers of $2$:

\begin{rmk}
  Notice in Proposition \ref{prop:2power} above that the leading diagonals of the
  matrix representations of $x^{2^n}$ and $y^{2^n}$ agree for all $n
  \ge 0$, since both elements are conjugate (see Lemma
  \ref{lem:crucform}(c)). Let $k = 2^n$. Recall that we can think of
  the elements in $H_k$ as matrices truncated at their $(k+1)$-st
  upper diagonal. Then the non-trivial group elements $x^k$ and $y^k$
  agree in $H_k$, since their leading diagonals coincide and are the
  $k$-th upper diagonals. (To separate these two elements in $H_k$,
  their first two non-trivial diagonals would have to survive under
  the truncation procedure.) This implies that
  \begin{equation} \label{eq:nontrivintersect}
  x_2 \Sigma(T_k) x_2^{-1} \cap \Sigma(T_k) \supset \{ x^k  \}. 
  \end{equation}
  Notice that condition (B) in the mixed Beauville structure implies
  the following property:
  \begin{itemize}
  \item[(B')] For all $g \in G\backslash H$: 
  $g \Sigma(T) g^{-1} \cap \Sigma(T) = \{ {\rm id} \}$,
  \end{itemize}
  since $\Sigma(T)$ is invariant under conjugation within $H$. But
  \eqref{eq:nontrivintersect} contradicts to (B') and we conclude that
  $(G_k,H_k,(x_0,x_1))$ cannot be a mixed Beauville structure if $k =
  2^n$.
\end{rmk}

To understand the first two non-trivial diagonals of {\em all powers} of $x$
and $y$ (not only the $2$-powers), we consider the binary presentation
of an arbitrary exponent $n \in {\mathbb N}$:
$$ n = 2^{k+j} \alpha_{k+j} + \cdots + 2^{k+1} \alpha_{k+1} + 2^k \alpha_k $$
with $\alpha_l \in \{ 0,1 \}$ for all $l$ and $\alpha_k=1$. Now define
\begin{equation} \label{eq:tn} 
t(n) = \begin{cases} 2 \alpha_1 + \alpha_0 & \text{if $k = 0$}, \\
  2^k \alpha_k & \text{if $k \ge 1$}.\end{cases} 
\end{equation}
Then $x^n$ is equal to $x^{t(n)}$ multiplied with certain higher
$2$-powers of $x$ (i.e., the powers $x^{2^{\alpha_l}}$ with $\alpha_l
=1$ and $l \ge \max\{2,k+1\}$). In view of Lemma
\ref{lem:crucform}(a), this multiplication does not change the first
two non-trivial diagonals of $x^{t(n)}$, which shows that the first
two non-trivial diagonals of $x^{t(n)}$ and $x^n$ agree. Using the
(easily computable) fact that
\begin{eqnarray*}
  x^3 &=& M_0([28,235,129],[46,138,451],\dots), \\
  y^3 &=& M_0([28,235,129],[9,90,377], \dots),
\end{eqnarray*}
this leads directly to the following result:

\begin{cor} \label{cor:first}
  The matrix representation of any power $x^n$ ($n \ge 1$) takes one
  of the following forms
  \begin{align*}
  & M_0([28,235,129],[29,211,263],\dots), 
  && M_0([28,235,129],[46,138,451],\dots) \\
  & M_{2^{\rm odd}-1}([51,89,196],[0,0,0],\dots),
  && M_{2^{{\rm even}+2}-1}([28,235,129],[0,0,0],\dots).
  \end{align*}
  The matrix representations of any power $y^n$ ($n \ge 1$)
  takes one of the following forms
  \begin{align*}
  & M_0([28,235,129],[58,3,445],\dots), 
  && M_0([28,235,129],[9,90,377],\dots) \\
  & M_{2^{\rm odd}-1}([51,89,196],[0,157,106],\dots),
  && M_{2^{{\rm even}+2}-1}([28,235,129],[39,208,186],\dots).
  \end{align*}
\end{cor}

Analogous results holds for the powers of the other four elements of
the two spherical systems of generators:

\begin{prop} \label{prop:second}
  The matrix representation of any power $x_0^n$ ($n \ge 1$) takes one
  of the following forms
  \begin{align*}
  & M_0([11,11,11],[17,17,17],\dots), 
  && M_0([11,11,11],[11,11,11],\dots) \\
  & M_{2^{\rm odd}-1}([26,26,26],[0,0,0],\dots),
  && M_{2^{{\rm even}+2}-1}([11,11,11],[0,0,0],\dots).
  \end{align*}
  The matrix representations of any power $y_0^n$ ($n \ge 1$)
  takes one of the following forms
  \begin{align*}
  & M_0([11,11,11],[44,219,177],\dots), 
  && M_0([11,11,11],[54,193,171],\dots) \\
  & M_{2^{\rm odd}-1}([26,26,26],[0,157,106],\dots),
  && M_{2^{{\rm even}+2}-1}([11,11,11],[61,202,160],\dots).
  \end{align*}
\end{prop}

\begin{prop} \label{prop:third}
  The matrix representation of any power $x_1^n$ ($n \ge 1$) takes one
  of the following forms
  \begin{align*}
  & M_0([23,224,138],[59,136,495],\dots), 
  && M_0([23,224,138],[28,88,341],\dots) \\
  & M_{2^{\rm odd}-1}([39,208,186],[0,0,0],\dots),
  && M_{2^{{\rm even}+2}-1}([23,224,138],[0,0,0],\dots).
  \end{align*}
  The matrix representations of any power $y_1^n$ ($n \ge 1$)
  takes one of the following forms
  \begin{align*}
  & M_0([23,224,138],[33,146,501],\dots), 
  && M_0([23,224,138],[6,66,335],\dots) \\
  & M_{2^{\rm odd}-1}([39,208,186],[0,106,247],\dots),
  && M_{2^{{\rm even}+2}-1}([23,224,138],[26,26,26],\dots).
  \end{align*}
\end{prop}

\section{Proof of property (C)}

The proof of property (C) for our triple $(G_k,H_k,T_k)$ is relatively
easy and follows solely from leading diagonal considerations. Since
every element in $H$ is a product of the elements
$x_0^{\pm 1},x_1^{\pm 1}$, we deduce first from Lemma \ref{lem:crucform}(a)
that the matrix representation of any element $x \in H$ takes one of
the following four forms: $M_0([0,0,0],\dots)$,
$M_0([11,11,11],\dots)$, $M_0([23,224,138],\dots)$ or
$M_0([28,235,129],\dots)$.

Using Lemma \ref{lem:crucform}, again, we obtain the following table:

\begin{center} 
\begin{tabular}{c|c|c}
$x \in H$ & $x \cdot x_2$ & $(x \cdot x_2)^2$ \\ \hline 
$M_0([0,0,0],\dots)$ & $M_0([46,68,217],\dots)$ & $M_1([41,67,222],\dots)$ 
\phantom{$\Big\Vert$}\\
$M_0([11,11,11],\dots)$ & $M_0([37,79,210],\dots)$ & $M_1([14,147,100],\dots])$
\phantom{$\Big\Vert$}\\
$M_0([23,224,138],\dots)$ & $M_0([57,164,83],\dots)$ & $M_1([20,137,126],\dots)$
\phantom{$\Big\Vert$}\\
$M_0([28,235,129],\dots)$ & $M_0([50,175,88],\dots)$ & $M_1([61,202,160],\dots)$
\phantom{$\Big\Vert$}
\end{tabular}
\end{center}

\smallskip

Now assume that $k \ge 2$. Since every element $g \in G_k \backslash
H_k$ is of the form $g = x x_2$ with $x \in H_k$, we conclude that
$g^2$ are truncations of matrices of one of the following four forms:
$M_1([41,67,222],\dots)$, $M_1([14,147,100],\dots])$,
$M_1([20,137,126],\dots)$ or $M_1([61,202,160],\dots)$. Notice that
the leading diagonal in the matrix of every such element $g^2$ is the
second upper diagonal.

On the other hand, since the leading diagonal of a matrix does not
change under conjugation (see Lemma \ref{lem:crucform}(c)), we
conclude from Corollary \ref{cor:first} and Propositions
\ref{prop:second} and \ref{prop:third} that the elements in
$\Sigma(T_k)$ are truncations of matrices of one of the following four
forms: $M_1([0,0,0],\dots)$, $M_1([51,89,196],\dots)$,
$M_1([26,26,26],\dots)$, $M_1([39,208,186],\dots)$.

Since these eight forms are all different, we conclude that $g^2
\not\in \Sigma(T_k)$ for all $g \in G_k \backslash H_k$. This shows
that property (C) in the definition of a mixed Beauville structure is
satisfied for all $k \ge 2$.

\section{Proof of property (B)}
\label{sec:propb}

In this section, we prove that our triples $(G_k,H_k,T_k)$ satisfy
property (B) of a mixed Beauville structure with the choice $g_0 =
x_2$, for all $k$ not a power of $2$. Recall that $x = (x_0 x_1)^{-1}$
and
$$ \Sigma(T) = \Sigma(x_0) \cup \Sigma(x_1) \cup \Sigma(x), $$
and
$$ x_2\Sigma(T)x_2^{-1} = \Sigma(y_0) \cup \Sigma(y_1) \cup \Sigma(y). $$
It follows immediately from inspection of the leading diagonals in
Corollary \ref{cor:first} and Propositions \ref{prop:second} and
\ref{prop:third} and the fact that these leading diagonals do not
change under conjugation (see Lemma \ref{lem:crucform}(c)) that, for
the pair $(x_0,y_1)$, we have
$$ \Sigma(x_0) \cap \Sigma(y_1) = \{ {\rm id} \}, $$
and that the same trivial intersection holds also for all other pairs
$(x_0,y)$, $(x_1,y_0)$, $(x_1,y)$, $(x,y_0)$ and $(x,y_1)$. So it only
remains to prove the trivial intersection
$$ \Sigma(x) \cap \Sigma(y) = \{ {\rm id} \}, $$
and analogous trivial intersection results for the pairs $(x_0,y_0)$
and $(x_1,y_1)$. For this, the consideration of the leading diagonal
is not sufficient and we have to study the behavior of the first two
non-trivial diagonals under conjugation. From now on, let $k$ be not a
power of $2$.

Note that 
$$ x,y = M_0(A_1 = [28,235,129],A_2,\dots) \mod G^k $$ 
with $A_2 = [29,211,263]$ or $A_2 = [58,3,445]$, respectively. Using
Lemma \ref{lem:crucform}(c), we see that $A_1$ does not change
under conjugation and that $A_2$ transforms under conjugation as follows:

\begin{equation*}
\begin{array}{c@{\hspace{1.8cm}}c@{\hspace{0.5cm}}c@{\hspace{1.8cm}}c}
  \Rnode{N1}{[29,211,263]} & \Rnode{N2}{[19,64,355]} 
& \Rnode{N5}{[58,3,445]} & \Rnode{N6}{[52,144,473]}\\[1.5cm]
  \Rnode{N3}{[19,64,355]} & \Rnode{N4}{[29,211,263]} 
& \Rnode{N7}{[52,144,473]} & \Rnode{N8}{[58,3,445]}
\end{array}
\psset{nodesep=0.3cm}
\everypsbox{\scriptstyle}
\ncLine{<->}{N1}{N2}\aput(0.5){{\rm Conj}(x_0^{\pm 1})}
\ncLine{<->}{N1}{N3}\Aput{{\rm Conj}(x_1^{\pm 1})}
\ncLine{<->}{N2}{N4}\Bput{{\rm Conj}(x_1^{\pm 1})}
\ncLine{<->}{N3}{N4}\aput(0.5){{\rm Conj}(x_0^{\pm 1})}
\ncLine{<->}{N5}{N6}\aput(0.5){{\rm Conj}(x_0^{\pm 1})}
\ncLine{<->}{N5}{N7}\Aput{{\rm Conj}(x_1^{\pm 1})}
\ncLine{<->}{N6}{N8}\Bput{{\rm Conj}(x_1^{\pm 1})}
\ncLine{<->}{N7}{N8}\aput(0.5){{\rm Conj}(x_0^{\pm 1})}
\end{equation*}

Since every element $h \in H_k$ is a product of the generators
$x_0^{\pm 1}, x_1^{\pm 1}$, we see that
$$
h x h^{-1} \neq h' y (h')^{-1} 
$$
for any pair $h,h' \in H_k$, since both elements differ in the second of
their first two non-trivial diagonals. Similarly, the conjugation scheme
for $A_2$ for the pair $x^3,y^3 = M_0(A_1=[28,235,129],A_2,\dots) \mod G^k$
reads as follows:

\begin{equation*}
\begin{array}{c@{\hspace{1.8cm}}c@{\hspace{0.5cm}}c@{\hspace{1.8cm}}c}
  \Rnode{N1}{[46,138,451]} & \Rnode{N2}{[32,25,423]} 
& \Rnode{N5}{[9,90,377]} & \Rnode{N6}{[7,201,285]}\\[1.5cm]
  \Rnode{N3}{[32,25,423]} & \Rnode{N4}{[46,138,451]} 
& \Rnode{N7}{[7,201,285]} & \Rnode{N8}{[9,90,377]}
\end{array}
\psset{nodesep=0.3cm}
\everypsbox{\scriptstyle}
\ncLine{<->}{N1}{N2}\aput(0.5){{\rm Conj}(x_0^{\pm 1})}
\ncLine{<->}{N1}{N3}\Aput{{\rm Conj}(x_1^{\pm 1})}
\ncLine{<->}{N2}{N4}\Bput{{\rm Conj}(x_1^{\pm 1})}
\ncLine{<->}{N3}{N4}\aput(0.5){{\rm Conj}(x_0^{\pm 1})}
\ncLine{<->}{N5}{N6}\aput(0.5){{\rm Conj}(x_0^{\pm 1})}
\ncLine{<->}{N5}{N7}\Aput{{\rm Conj}(x_1^{\pm 1})}
\ncLine{<->}{N6}{N8}\Bput{{\rm Conj}(x_1^{\pm 1})}
\ncLine{<->}{N7}{N8}\aput(0.5){{\rm Conj}(x_0^{\pm 1})}
\end{equation*}

Comparison of the second non-trivial diagonals shows that we have
\begin{equation} \label{eq:xnneqxm}
h x^n h^{-1} \neq h' y^m (h')^{-1} 
\end{equation}
for any pair $h,h' \in H_k$ and every $n,m$ with
$t(n),t(m) \in \{1,3\}$, where $t(n)$ was defined in \eqref{eq:tn}.

Next, let us look at the conjugation scheme for $A_2$ for any pair
$x^{2^r},y^{2^r} = M_{2^r-1}(A_1=[51,89,196],A_2) \mod G^k$, where $r
\ge 1$ is odd:

\begin{equation*}
\begin{array}{c@{\hspace{1.8cm}}c@{\hspace{0.5cm}}c@{\hspace{1.8cm}}c}
  \Rnode{N1}{[0,0,0]} & \Rnode{N2}{[0,247,157]} 
& \Rnode{N5}{[0,157,106]} & \Rnode{N6}{[0,106,247]}\\[1.5cm]
  \Rnode{N3}{[0,247,157]} & \Rnode{N4}{[0,0,0]} 
& \Rnode{N7}{[0,106,247]} & \Rnode{N8}{[0,157,106]}
\end{array}
\psset{nodesep=0.3cm}
\everypsbox{\scriptstyle}
\ncLine{<->}{N1}{N2}\aput(0.5){{\rm Conj}(x_0^{\pm 1})}
\ncLine{<->}{N1}{N3}\Aput{{\rm Conj}(x_1^{\pm 1})}
\ncLine{<->}{N2}{N4}\Bput{{\rm Conj}(x_1^{\pm 1})}
\ncLine{<->}{N3}{N4}\aput(0.5){{\rm Conj}(x_0^{\pm 1})}
\ncLine{<->}{N5}{N6}\aput(0.5){{\rm Conj}(x_0^{\pm 1})}
\ncLine{<->}{N5}{N7}\Aput{{\rm Conj}(x_1^{\pm 1})}
\ncLine{<->}{N6}{N8}\Bput{{\rm Conj}(x_1^{\pm 1})}
\ncLine{<->}{N7}{N8}\aput(0.5){{\rm Conj}(x_0^{\pm 1})}
\end{equation*}

Again, this shows that we have \eqref{eq:xnneqxm} for any pair $h,h'
\in H_k$ and every $n,m$ with $t(n)=t(m)=2^r < k$ and odd $r \ge 1$.
Moreover, \eqref{eq:xnneqxm} also holds for any choice of $n,m$ such
that
\begin{itemize}
\item[(i)] one of $t(n),t(m)$ is in $\{1,3\}$ and the other is of the
  form $2^r$ with odd $r \ge 1$, or
\item[(ii)] $t(n) = 2^{r_1} < k$ and $t(m) = 2^{r_2} < k$ with $r_1,
  r_2 \ge 1$ both odd and $r_1 \neq r_2$,
\end{itemize}
since then the number of first upper vanishing diagonals of $h x^n
h^{-1}$ and $h' y^m (h')^{-1}$ do not agree.

Finally, we have the following conjugation scheme for $A_2$ for any pair
$x^{2^r},y^{2^r}=(A_1=[28,235,129],A_2) \mod G^k$ with even $r \ge 2$:

\begin{equation*}
\begin{array}{c@{\hspace{1.8cm}}c@{\hspace{0.5cm}}c@{\hspace{1.8cm}}c}
  \Rnode{No1}{[0,0,0]} & \Rnode{No2}{[14,147,100]} 
& \Rnode{No5}{[39,208,186]} & \Rnode{No6}{[41,67,222]}\\[1.5cm]
  \Rnode{No3}{[14,147,100]} & \Rnode{No4}{[0,0,0]} 
& \Rnode{No7}{[41,67,222]} & \Rnode{No8}{[39,208,186]}
\end{array}
\psset{nodesep=0.3cm}
\everypsbox{\scriptstyle}
\ncLine{<->}{No1}{No2}\aput(0.5){{\rm Conj}(x_0^{\pm 1})}
\ncLine{<->}{No1}{No3}\Aput{{\rm Conj}(x_1^{\pm 1})}
\ncLine{<->}{No2}{No4}\Bput{{\rm Conj}(x_1^{\pm 1})}
\ncLine{<->}{No3}{No4}\aput(0.5){{\rm Conj}(x_0^{\pm 1})}
\ncLine{<->}{No5}{No6}\aput(0.5){{\rm Conj}(x_0^{\pm 1})}
\ncLine{<->}{No5}{No7}\Aput{{\rm Conj}(x_1^{\pm 1})}
\ncLine{<->}{No6}{No8}\Bput{{\rm Conj}(x_1^{\pm 1})}
\ncLine{<->}{No7}{No8}\aput(0.5){{\rm Conj}(x_0^{\pm 1})}
\end{equation*}

Combining all above results shows that we have \eqref{eq:xnneqxm} for all
$n,m \ge 1$ with $t(n),t(m) \le k$. Note that $x^n = y^n = {\rm id}$ for all
$n \ge 1$ with $t(n) > k$, so we conclude that
$$ \Sigma(x) \cap \Sigma(y) = \{ {\rm id} \}. $$  
The corresponding commutator schemes for the pairs $(x_0,y_0)$ and
$(x_1,y_1)$ are listed in Appendices \ref{app:x0y0} and
\ref{app:x1y1} below, finishing the proof of
$$ x_2 \Sigma(T) x_2^{-1} \cap \Sigma(T) = \{ {\rm id} \}. $$ 

\section{Bringing everything together}

Is is obvious that our triples $u_k = (G_k,H_k,T_k)$ satisfy property (A) of
a mixed Beauville structure. Since the previous two sections show the
validity of properties (B) and (C) if $k$ is not a power of $2$, we conclude
that these triples are mixed Beauville structures. 

\smallskip

Next, we use the following fact: {\em Assume that $\Gamma$ is a finite group
with finite presentation, i.e.,
$$ \Gamma = \langle g_0,\dots,g_k \mid r_1(g_0,\dots,g_k) = {\rm id},\dots,
r_l(g_0,\dots,g_k) = {\rm id} \rangle. $$
For $0 \le i \le k$, let
$$ g_i' = w_i(g_0,\dots,g_k) $$
and
$$ r_i'(g_0,\dots,g_k) = r_i(w_0(g_0,\dots,g_k),\dots,w_k(g_0,\dots,g_k)). $$
Let $\Gamma'$ be the group defined by
\begin{multline*}
\Gamma' = \langle g_0,\dots,g_k \mid r_1(g_0,\dots,g_k)
= r_1'(g_0,\dots,g_k) = {\rm id},\dots,\\
r_l(g_0,\dots,g_k) = r_l'(g_0,\dots,g_k) = {\rm id} \rangle. 
\end{multline*}
Then there exists a unique homomorphism
$\psi: \Gamma \to \Gamma$ with $\psi(g_i) = g_i'$ ($0 \le i \le k$) if
and only if $|\Gamma| = |\Gamma'|$.}

\smallskip

In view of \cite[p. 2782]{PV}, it is easily checked that $G$ is
canonically isomorphic to
$$ \langle x_0,x_1,x_2 \mid r_1(x_0,x_1,x_2) = r_2(x_0,x_1,x_2) = 
r_3(x_0,x_1,x_2) = {\rm id} \rangle, $$
with
\begin{eqnarray*}
  r_1(x_0,x_1,x_2) &=& x_2 x_1 x_2 x_0 x_1 x_0, \\
  r_2(x_0,x_1,x_2) &=& x_2 x_0^{-1} x_2 x_1^{-1} x_0^{-1} x_1, \\
  r_3(x_0,x_1,x_2) &=& x_2^2 x_1^{-1} x_0^{-1} x_1^{-1} x_0,
\end{eqnarray*}
and that the quotient $G_3$ is canonically isomorphic to
\begin{multline*}
\langle x_0,x_1,x_2 \mid r_1(x_0,x_1,x_2) = r_2(x_0,x_1,x_2) = 
r_3(x_0,x_1,x_2) = {\rm id}, \\
[x_1,x_0,x_0,x_0] = [x_1,x_0,x_0,x_1] = [x_1,x_0,x_0,x_2] = {\rm id}
\rangle.
\end{multline*}
Using the above criterion, a straightforward MAGMA calculation shows
that there exists a unique automorphism $\psi: G_3 \to G_3$ with
$\psi(x_0) = x_0^{-1}$, $\psi(x_1) = x_1^{-1}$ and $\psi(x_2) =
x_0^{-1} x_2 x_0$. This shows that $S(u_3)$ is a {\em real} Beauville
surface of mixed type.

\smallskip

Finally, recall from \cite[p. 2781]{PV} that $H$ is canonically isomorphic
to 
$$ \langle x_0,x_1 \mid r_3(x_0,x_1) = r_4(x_0,x_1) = 
r_5(x_0,x_1) = {\rm id} \rangle, $$
with
\begin{eqnarray*}
  r_3(x_0,x_1) &=& (x_1 x_0)^3 x_1^{-3} x_0^{-3}, \\
  r_4(x_0,x_1) &=& x_1 x_0^{-1} x_1^{-1} x_0^{-3} x_1^2 x_0^{-1} x_1 x_0 x_1, \\
  r_5(x_0,x_1) &=& x_1^3 x_0^{-1} x_1 x_0 x_1 x_0^2 x_1^2 x_0 x_1 x_0,
\end{eqnarray*}
MAGMA calculations show that for any choice 
\begin{multline*}
(y_0,y_1) \in \{ (x_0^{-1},x_1^{-1}), (x_1 x_0,x_0^{-1}),
(x_1^{-1},x_1 x_0), \\ (x_1^{-1},x_0^{-1}), (x_0^{-1},x_1 x_0),
(x_1 x_0,x_1^{-1}) \},
\end{multline*}
we have
\begin{multline*}
| \langle x_0,x_1 \mid r_3(x_0,x_1) = r_3(y_0,y_1) = {\rm id},
r_4(x_0,x_1) = r_4(y_0,y_1) = {\rm id}, \\
r_5(x_0,x_1) = r_5(y_0,y_1) = {\rm id} \rangle | = 3072.
\end{multline*}
Since we have $|H_k| \ge 8192$ for $k > 3$ not a power of $2$, there
cannot be a homomorphism $\psi: H_k \to H_k$ with $\psi(x_0) = y_0$
and $\psi(x_1) = y_1$ by the above criterion, showing that $S(u_k)$
cannot be biholomorphic to $\overline{S(u_k)}$. This finishes the proof
of Theorem \ref{thm:main}. \qed

\begin{rmk} \label{rmk:reflection} The fact that the triples
  $(G_k,H_k,T_k)$ are mixed Beauville structures (for $k$ not a power
  of $2$) is very remarkable. Let us reflect -- by looking back at the
  proof of property (B) -- why this result is so surprising:

  We know that for indices up to order $k \le 100$ we have
  $$ |G_{k+1}| = \begin{cases} 8 |G_k|, & \text{if $k \equiv 0,1 \mod 3$,}\\ 
    4 |G_k|, & \text{if $k \equiv 2 \mod 3$,} \end{cases} $$ which
  gives strong evidence that this should hold for all indices $k \in
  {\mathbb N}$ (see the {\em finite width $3$ conjecture} in
  \cite[Conjecture 1]{PV}).

  This means that for any $k \le 99$, $k \equiv 2 \mod 3$ and $A_1$,
  there are {\em at most four} different choices $A_2$ such that
  $(A_1,A_2)$ represent the first two non-trivial diagonals of matrix
  representations $M_k(A_1,A_2,\dots)$ of elements in $G$. On the
  other hand, it follows from the arguments in the proof of property
  (B) that we need {\em at least four} such possibilities to guarantee
  that $\Sigma(x) \cap \Sigma(y) = \{ {\rm id} \}$ (and to derive
  analogous results for the pairs $(x_0,y_0)$ and
  $(x_1,y_1)$). Moreover, these four possibilities must appear in the
  right combinations in all the conjugation schemes to guarantee the
  required trivial intersections.

  Moreover, for any given $A_1$, we have {\em at most eight} choices
  $A_2$ such that $M_0(A_1,A_2,\dots)$ are matrix representations of
  elements in $G$. On the other hand, our considerations in the
  previous section show that we need {\em at least eight} such choices
  to guarantee that
  $$ \{ h x h^{-1}, h x^3 h^{-1}\} \cap \{ h' y h^{-1}, h y^3 (h')^{-1} \} =
  \{ {\rm id} \}, $$
  for all choices of $h,h' \in H$. 

  This shows that the conjectured {\em finite width $3$} property of the
  infinite group $G$ implies a very tight situation, which leaves
  ''just enough room'' to allow the mixed Beauville structures (for $k$
  not a power of $2$). 
\end{rmk}

\begin{appendix}

\section{Conjugation schemes for the pairs $x_0^n,y_0^n$}
\label{app:x0y0}

The notation in the conjugation schemes is the same as in Section
\ref{sec:propb}. The results in this appendix show that
$$ \Sigma(x_0) \cap \Sigma(y_0) = \{ {\rm id} \}. $$

\begin{itemize}
\item[(a)] For the pair $x_0,y_0 = M_0(A_1=[11,11,11],A_2,\dots) \mod
  G^k$:

  \begin{equation*}
    \begin{array}{c@{\hspace{1.8cm}}c@{\hspace{0.5cm}}c@{\hspace{1.8cm}}c}
      \Rnode{N1}{[17,17,17]} & \Rnode{N2}{[17,17,17]} 
      & \Rnode{N5}{[44,219,177]} & \Rnode{N6}{[44,219,177]}\\[1.5cm]
      \Rnode{N3}{[31,130,117]} & \Rnode{N4}{[31,130,117]} 
      & \Rnode{N7}{[34,72,213]} & \Rnode{N8}{[34,72,213]}
    \end{array}
    \psset{nodesep=0.3cm}
    \everypsbox{\scriptstyle}
    \ncLine{<->}{N1}{N2}\aput(0.5){{\rm Conj}(x_0^{\pm 1})}
    \ncLine{<->}{N1}{N3}\Aput{{\rm Conj}(x_1^{\pm 1})}
    \ncLine{<->}{N2}{N4}\Bput{{\rm Conj}(x_1^{\pm 1})}
    \ncLine{<->}{N3}{N4}\aput(0.5){{\rm Conj}(x_0^{\pm 1})}
    \ncLine{<->}{N5}{N6}\aput(0.5){{\rm Conj}(x_0^{\pm 1})}
    \ncLine{<->}{N5}{N7}\Aput{{\rm Conj}(x_1^{\pm 1})}
    \ncLine{<->}{N6}{N8}\Bput{{\rm Conj}(x_1^{\pm 1})}
    \ncLine{<->}{N7}{N8}\aput(0.5){{\rm Conj}(x_0^{\pm 1})}
  \end{equation*}
\item[(b)] For the pair $x_0^3,y_0^3 = M_0(A_1=[11,11,11],A_2,\dots) \mod
  G^k$:

  \begin{equation*}
    \begin{array}{c@{\hspace{1.8cm}}c@{\hspace{0.5cm}}c@{\hspace{1.8cm}}c}
      \Rnode{N1}{[11,11,11]} & \Rnode{N2}{[11,11,11]} 
      & \Rnode{N5}{[54,193,171]} & \Rnode{N6}{[54,193,171]}\\[1.5cm]
      \Rnode{N3}{[5,152,111]} & \Rnode{N4}{[5,152,111]} 
      & \Rnode{N7}{[56,82,207]} & \Rnode{N8}{[56,82,207]}
    \end{array}
    \psset{nodesep=0.3cm}
    \everypsbox{\scriptstyle}
    \ncLine{<->}{N1}{N2}\aput(0.5){{\rm Conj}(x_0^{\pm 1})}
    \ncLine{<->}{N1}{N3}\Aput{{\rm Conj}(x_1^{\pm 1})}
    \ncLine{<->}{N2}{N4}\Bput{{\rm Conj}(x_1^{\pm 1})}
    \ncLine{<->}{N3}{N4}\aput(0.5){{\rm Conj}(x_0^{\pm 1})}
    \ncLine{<->}{N5}{N6}\aput(0.5){{\rm Conj}(x_0^{\pm 1})}
    \ncLine{<->}{N5}{N7}\Aput{{\rm Conj}(x_1^{\pm 1})}
    \ncLine{<->}{N6}{N8}\Bput{{\rm Conj}(x_1^{\pm 1})}
    \ncLine{<->}{N7}{N8}\aput(0.5){{\rm Conj}(x_0^{\pm 1})}
  \end{equation*}
\item[(c)] For the pair $x_0^{2^r},y_0^{2^r} =
  M_{2^r-1}(A_1=[26,26,26],A_2,\dots) \mod G^k$, where $r \ge 1$ is
  odd:

  \begin{equation*}
    \begin{array}{c@{\hspace{1.8cm}}c@{\hspace{0.5cm}}c@{\hspace{1.8cm}}c}
      \Rnode{N1}{[0,0,0]} & \Rnode{N2}{[0,0,0]} 
      & \Rnode{N5}{[0,157,106]} & \Rnode{N6}{[0,157,106]}\\[1.5cm]
      \Rnode{N3}{[0,106,247]} & \Rnode{N4}{[0,106,247]} 
      & \Rnode{N7}{[0,247,157]} & \Rnode{N8}{[0,247,157]}
    \end{array}
    \psset{nodesep=0.3cm}
    \everypsbox{\scriptstyle}
    \ncLine{<->}{N1}{N2}\aput(0.5){{\rm Conj}(x_0^{\pm 1})}
    \ncLine{<->}{N1}{N3}\Aput{{\rm Conj}(x_1^{\pm 1})}
    \ncLine{<->}{N2}{N4}\Bput{{\rm Conj}(x_1^{\pm 1})}
    \ncLine{<->}{N3}{N4}\aput(0.5){{\rm Conj}(x_0^{\pm 1})}
    \ncLine{<->}{N5}{N6}\aput(0.5){{\rm Conj}(x_0^{\pm 1})}
    \ncLine{<->}{N5}{N7}\Aput{{\rm Conj}(x_1^{\pm 1})}
    \ncLine{<->}{N6}{N8}\Bput{{\rm Conj}(x_1^{\pm 1})}
    \ncLine{<->}{N7}{N8}\aput(0.5){{\rm Conj}(x_0^{\pm 1})}
  \end{equation*} 
\item[(d)] For the pair $x_0^{2^r},y_0^{2^r} =
  M_{2^r-1}(A_1=[11,11,11],A_2,\dots) \mod G^k$, where $r \ge 2$ is
  even:
 
  \begin{equation*}
    \begin{array}{c@{\hspace{1.8cm}}c@{\hspace{0.5cm}}c@{\hspace{1.8cm}}c}
      \Rnode{N1}{[0,0,0]} & \Rnode{N2}{[0,0,0]} 
      & \Rnode{N5}{[61,202,160]} & \Rnode{N6}{[61,202,160]}\\[1.5cm]
      \Rnode{N3}{[14,147,100]} & \Rnode{N4}{[14,147,100]} 
      & \Rnode{N7}{[51,89,196]} & \Rnode{N8}{[51,89,196]}
    \end{array}
    \psset{nodesep=0.3cm}
    \everypsbox{\scriptstyle}
    \ncLine{<->}{N1}{N2}\aput(0.5){{\rm Conj}(x_0^{\pm 1})}
    \ncLine{<->}{N1}{N3}\Aput{{\rm Conj}(x_1^{\pm 1})}
    \ncLine{<->}{N2}{N4}\Bput{{\rm Conj}(x_1^{\pm 1})}
    \ncLine{<->}{N3}{N4}\aput(0.5){{\rm Conj}(x_0^{\pm 1})}
    \ncLine{<->}{N5}{N6}\aput(0.5){{\rm Conj}(x_0^{\pm 1})}
    \ncLine{<->}{N5}{N7}\Aput{{\rm Conj}(x_1^{\pm 1})}
    \ncLine{<->}{N6}{N8}\Bput{{\rm Conj}(x_1^{\pm 1})}
    \ncLine{<->}{N7}{N8}\aput(0.5){{\rm Conj}(x_0^{\pm 1})}
  \end{equation*}
\end{itemize}

\section{Conjugation schemes for the pairs $x_1^n,y_1^n$}
\label{app:x1y1}

The notation in the conjugation schemes is the same as in Section
\ref{sec:propb}. The results in this appendix show that
$$ \Sigma(x_1) \cap \Sigma(y_1) = \{ {\rm id} \}. $$

\begin{itemize}
\item[(a)] For the pair $x_1,y_1 = M_0(A_1=[23,224,138],A_2,\dots) \mod
  G^k$:

  \begin{equation*}
    \begin{array}{c@{\hspace{1.8cm}}c@{\hspace{0.5cm}}c@{\hspace{1.8cm}}c}
      \Rnode{N1}{[59,136,495]} & \Rnode{N2}{[53,27,395]} 
      & \Rnode{N5}{[33,146,501]} & \Rnode{N6}{[47,1,401]}\\[1.5cm]
      \Rnode{N3}{[59,136,495]} & \Rnode{N4}{[53,27,395]} 
      & \Rnode{N7}{[33,146,501]} & \Rnode{N8}{[47,1,401]}
    \end{array}
    \psset{nodesep=0.3cm}
    \everypsbox{\scriptstyle}
    \ncLine{<->}{N1}{N2}\aput(0.5){{\rm Conj}(x_0^{\pm 1})}
    \ncLine{<->}{N1}{N3}\Aput{{\rm Conj}(x_1^{\pm 1})}
    \ncLine{<->}{N2}{N4}\Bput{{\rm Conj}(x_1^{\pm 1})}
    \ncLine{<->}{N3}{N4}\aput(0.5){{\rm Conj}(x_0^{\pm 1})}
    \ncLine{<->}{N5}{N6}\aput(0.5){{\rm Conj}(x_0^{\pm 1})}
    \ncLine{<->}{N5}{N7}\Aput{{\rm Conj}(x_1^{\pm 1})}
    \ncLine{<->}{N6}{N8}\Bput{{\rm Conj}(x_1^{\pm 1})}
    \ncLine{<->}{N7}{N8}\aput(0.5){{\rm Conj}(x_0^{\pm 1})}
  \end{equation*}
\item[(b)] For the pair $x_1^3,y_1^3 = M_0(A_1=[23,224,138],A_2,\dots) \mod
  G^k$:

  \begin{equation*}
    \begin{array}{c@{\hspace{1.8cm}}c@{\hspace{0.5cm}}c@{\hspace{1.8cm}}c}
      \Rnode{N1}{[28,88,341]} & \Rnode{N2}{[18,203,305]} 
      & \Rnode{N5}{[6,66,335]} & \Rnode{N6}{[8,209,299]}\\[1.5cm]
      \Rnode{N3}{[28,88,341]} & \Rnode{N4}{[18,203,305]} 
      & \Rnode{N7}{[6,66,335]} & \Rnode{N8}{[8,209,299]}
    \end{array}
    \psset{nodesep=0.3cm}
    \everypsbox{\scriptstyle}
    \ncLine{<->}{N1}{N2}\aput(0.5){{\rm Conj}(x_0^{\pm 1})}
    \ncLine{<->}{N1}{N3}\Aput{{\rm Conj}(x_1^{\pm 1})}
    \ncLine{<->}{N2}{N4}\Bput{{\rm Conj}(x_1^{\pm 1})}
    \ncLine{<->}{N3}{N4}\aput(0.5){{\rm Conj}(x_0^{\pm 1})}
    \ncLine{<->}{N5}{N6}\aput(0.5){{\rm Conj}(x_0^{\pm 1})}
    \ncLine{<->}{N5}{N7}\Aput{{\rm Conj}(x_1^{\pm 1})}
    \ncLine{<->}{N6}{N8}\Bput{{\rm Conj}(x_1^{\pm 1})}
    \ncLine{<->}{N7}{N8}\aput(0.5){{\rm Conj}(x_0^{\pm 1})}
  \end{equation*}
\item[(c)] For the pair $x_1^{2^r},y_1^{2^r} =
  M_{2^r-1}(A_1=[39,208,186],A_2,\dots) \mod G^k$, where $r \ge 1$ is
  odd:

  \begin{equation*}
    \begin{array}{c@{\hspace{1.8cm}}c@{\hspace{0.5cm}}c@{\hspace{1.8cm}}c}
      \Rnode{N1}{[0,0,0]} & \Rnode{N2}{[0,157,106]} 
      & \Rnode{N5}{[0,106,247]} & \Rnode{N6}{[0,247,157]}\\[1.5cm]
      \Rnode{N3}{[0,0,0]} & \Rnode{N4}{[0,157,106]} 
      & \Rnode{N7}{[0,106,247]} & \Rnode{N8}{[0,247,157]}
    \end{array}
    \psset{nodesep=0.3cm}
    \everypsbox{\scriptstyle}
    \ncLine{<->}{N1}{N2}\aput(0.5){{\rm Conj}(x_0^{\pm 1})}
    \ncLine{<->}{N1}{N3}\Aput{{\rm Conj}(x_1^{\pm 1})}
    \ncLine{<->}{N2}{N4}\Bput{{\rm Conj}(x_1^{\pm 1})}
    \ncLine{<->}{N3}{N4}\aput(0.5){{\rm Conj}(x_0^{\pm 1})}
    \ncLine{<->}{N5}{N6}\aput(0.5){{\rm Conj}(x_0^{\pm 1})}
    \ncLine{<->}{N5}{N7}\Aput{{\rm Conj}(x_1^{\pm 1})}
    \ncLine{<->}{N6}{N8}\Bput{{\rm Conj}(x_1^{\pm 1})}
    \ncLine{<->}{N7}{N8}\aput(0.5){{\rm Conj}(x_0^{\pm 1})}
  \end{equation*} 
\item[(d)] For the pair $x_1^{2^r},y_1^{2^r} =
  M_{2^r-1}(A_1=[23,224,138],A_2,\dots) \mod G^k$, where $r \ge 2$ is
  even:
 
  \begin{equation*}
    \begin{array}{c@{\hspace{1.8cm}}c@{\hspace{0.5cm}}c@{\hspace{1.8cm}}c}
      \Rnode{N1}{[0,0,0]} & \Rnode{N2}{[14,147,100]} 
      & \Rnode{N5}{[26,26,26]} & \Rnode{N6}{[20,137,126]}\\[1.5cm]
      \Rnode{N3}{[0,0,0]} & \Rnode{N4}{[14,147,100]} 
      & \Rnode{N7}{[26,26,26]} & \Rnode{N8}{[20,137,126]}
    \end{array}
    \psset{nodesep=0.3cm}
    \everypsbox{\scriptstyle}
    \ncLine{<->}{N1}{N2}\aput(0.5){{\rm Conj}(x_0^{\pm 1})}
    \ncLine{<->}{N1}{N3}\Aput{{\rm Conj}(x_1^{\pm 1})}
    \ncLine{<->}{N2}{N4}\Bput{{\rm Conj}(x_1^{\pm 1})}
    \ncLine{<->}{N3}{N4}\aput(0.5){{\rm Conj}(x_0^{\pm 1})}
    \ncLine{<->}{N5}{N6}\aput(0.5){{\rm Conj}(x_0^{\pm 1})}
    \ncLine{<->}{N5}{N7}\Aput{{\rm Conj}(x_1^{\pm 1})}
    \ncLine{<->}{N6}{N8}\Bput{{\rm Conj}(x_1^{\pm 1})}
    \ncLine{<->}{N7}{N8}\aput(0.5){{\rm Conj}(x_0^{\pm 1})}
  \end{equation*}
\end{itemize}

\end{appendix}

\end{document}